\documentclass{amsart}   
\usepackage{amsfonts}
\usepackage{amsmath,amssymb}
\usepackage{epsfig}
\newcommand{\MSOL}{\mathbf{MSOL}}
\newif\ifskip
\skipfalse
\newif\ifcomments
\commentsfalse
\newif\ifcpc
\cpcfalse 
\newtheorem{theorem}{Theorem}
\newtheorem{proposition}[theorem]{Proposition}
\newtheorem{remark}[theorem]{Remark}
\newtheorem{corollary}[theorem]{Corollary}
\newtheorem{problem}[theorem]{Problem}
\ifcpc

\newenvironment{proof}[1][Proof]{\noindent\textbf{#1.} }{\ \rule{0.5em}{0.5em}}
\else
\fi 
\begin{document}
\begin{titlepage}
\title[Enumeration of Vertex Induced Subgraphs]
{The Enumeration of Vertex Induced Subgraphs with respect to the
Number of Components}
\author[P. Tittmann]{P. Tittmann}
\author[I. Averbouch]{I. Averbouch$^{*}$}
\author[J.A. Makowsky]{J.A. Makowsky$^{**}$}

\thanks{$^{*}$
Partially supported by a grant of the Graduate School of the Technion--Israel Institute of Technology}
\thanks{$^{**}$
Partially supported by a grant of the Fund for
Promotion of Research of the Technion--Israel Institute of
Technology and grant ISF 1392/07 of the Israel Science Foundation (2007-2010)}
\email[P. Tittmann]{peter@HTWM.De}
\email[I. Averbouch]{ailia@cs.technion.ac.il}
\email[J.A. Makowsky]{janos@cs.technion.ac.il}
\date{January 20, 2009}

\address[P. Tittmann]{
Fachbereich Mathematik--Physik--Informatik,
\newline
Hochschule Mittweida, 
Mittweida, Germany 
}
\address[I. Averbouch and J.A. Makowsky]{
Department of Computer Science,
\newline
Technion--Israel Institute of Technology,
32000 Haifa, Israel
}

\begin{abstract}
Inspired by the study of community structure in connection
networks,
we introduce the graph polynomial
$Q\left( G;x,y\right)$, 
the bivariate generating function
which counts the number of connected components
in induced subgraphs.

We give a recursive definition of 
$Q\left( G;x,y\right)$ 
using vertex deletion, vertex contraction and
deletion of a vertex together with its neighborhood
and prove a universality property.
We relate 
$Q\left( G;x,y\right)$ 
to other known graph invariants and graph polynomials, among them 
partition functions, the Tutte polynomial, the independence and matching polynomials,
and the universal edge elimination polynomial introduced by I. Averbouch, B. Godlin
and J.A. Makowsky (2008).

We show that
$Q\left( G;x,y\right)$  is vertex reconstructible in the sense of Kelly and Ulam,
discuss its use in computing residual connectedness reliability.
Finally we show that 
the computation of $Q\left( G;x,y\right)$  
is $\sharp \mathbf{P}$-hard,
but  Fixed Parameter Tractable for graphs of bounded tree-width and clique-width.

\end{abstract}
\end{titlepage}
\maketitle

\tiny
\tableofcontents
\newpage
\normalsize

\section{Introduction}

\subsection{Motivation: Community Structure in Networks}
In the last decade
stochastic social networks have been analyzed mathematically from various
points of view. Understanding such networks sheds light on many questions
arising in biology, epidemology, sociology and large computer networks.
Researchers have concentrated particularly on a few properties 
that seem to be common to many networks: the small-world property, 
power-law degree distributions, 
and network transitivity, 
For  a broad view on the structure and dynamics of networks,
see \cite{bk:NewmanBarabasiWatts06}.
M. Girvan and  M.E.J. Newman, \cite{ar:GirvanNewman02},
highlight another property that is found in many networks, 
the property of {\em community structure}, 
in which network nodes are joined together in tightly knit groups, 
between which there are only looser connections.

Motivated by 
\cite{ar:Newman04a}, and
the first author's involvement in a project studying social networks,
we were led to study the graph parameter 
$q_{ij}\left( G\right) $, the number of vertex subsets $X\subseteq V$ with 
$i$ vertices such that $G\left[ X\right] $ has exactly $j$ components.
$q_{ij}\left( G\right) $, counts the number of degenerated communities which consist of $i$
members, and which split into $j$  isolated subcommunities.

The ordinary bivariate generating function associated with 
$q_{ij}\left( G\right) $
is the two-variable
graph polynomial
\[
Q\left( G;x,y\right) =\sum_{i=0}^{n}\sum_{j=0}^{n}q_{ij}\left( G\right)
x^{i}y^{j}. 
\]
We call $Q\left( G;x,y\right) $ the 
\emph{subgraph component polynomial} of $G$. 
The coefficient of $y^{k}$ in $Q\left( G;x,y\right) $ is the ordinary
generating function for the number of vertex sets that induce a subgraph of 
$G$ with exactly $k$ components.

\subsection{$Q(G;x,y)$ as a Graph Polynomial}
There is an abundance of graph polynomials studied in the literature,
and slowly there is a framework emerging, 
\cite{pr:Makowsky06,ar:MakowskyZoo,pr:GodlinKatzMakowsky07},
which allows to compare graph polynomials with respect to their
ability to distinguish graphs, to encode other graph polynomials or
numeric graph invariants, and their computational complexity.
In this paper we study
the {\emph subgraph component polynomial}
$Q\left( G;x,y\right)$ as a graph polynomial in its own right
and explore its properties within this emerging framework.

Like the bivariate Tutte polynomial,
see \cite[Chapter 10]{bk:Bollobas99},
the  polynomial
$Q\left( G;x,y\right)$  has several remarkable properties.
However, its distinguishing power is quite different from the
Tutte polynomial and other well studied polynomials.

Our main findings are:
\begin{itemize}
\item
$Q\left( G;x,y\right)$ distinguishes graphs which cannot be distinguished
by the matching polynomial, the Tutte polynomial, the characteristic
polynomial, or the bivariate chromatic polynomial introduced in
\cite{ar:DPT03} (Section \ref{se:distinct}).
\item
Nevertheless, we construct an infinite family of pairs of graphs
which cannot be pairwise distinguished by $Q\left( G;x,y\right)$ 
(Proposition \ref{prop:inftrees}). 
\item
The Tutte polynomial, 
satisfies a linear recurrence relation
with respect to edge deletion and edge contraction, and is universal
in this respect.
$Q\left( G;x,y\right)$  
also 
satisfies a linear recurrence relation, but with respect to three
kinds of vertex elimination, and is universal in this respect.
(Theorems \ref{theo_decom} and \ref{th:universal}).
\item
A graph polynomial in three indeterminates, 
$\xi(G;x,y,z)$,  
which satisfies a linear recurrence relation
with respect to three kinds of edge elimation, and which is universal
in this respect, was introduced in 
\cite{pr:AverbouchGodlinMakowsky07,ar:AverbouchGodlinMakowsky08}.
It subsumes both the Tutte polynomial and the matching polynomial.
For a line graph $L(G)$ of a graph $G$, we have
$Q\left( L(G);x,y\right)$   is a substitution instance of
$\xi(G;x,y,z)$
(Theorem \ref{thm:q_xi}).
\item
For fixed positive integer $n$ the univariate polynomial
$Q(G;x,n)$ can be interpreted as counting weighted homomorphisms,
\cite{ar:DyerGreenhill2000},
and is related to the Widom-Rowlinson model for $n$ particles 
(Theorem \ref{thm:hom}).
\item
$Q(G;x,y)$ is reconstructible from its vertex deletion deck
in the sense of \cite{ar:BondyHemminger77,ar:Bondy91} 
(Theorem \ref{th:reconstruct}).
\item
$Q(G;x,y)$ can be used (Section \ref{se:random}), to compute the
probability $P_{k}\left( G\right) $ that a vertex induced subgraph of $G$
has exactly $k$ components from the subgraph polynomial.
For $k=1$ this is known as the
{\em residual connectedness reliability}
(Section \ref{se:random}).
\item
Also like for the Tutte polynomial, cf. \cite{ar:JaegerVertiganWelsh90},
$Q(G;x_0,y_0)$  has the {\em Difficult Point Property}, i.e. it
is $\sharp \mathbf{P}$-hard to compute for 
all fixed values of $(x_0,y_0) \in \mathbb{R}^2 -E$ 
where $E$ is a semi-algebraic set of lower dimension
(Theorem \ref{th:dpp}).
In \cite{ar:MakowskyZoo} it is conjectured that
the Difficult Point Property holds for a wide class of graph polynomials,
the graph polynomials definable in Monadic Second Order Logic.
The conjecture has been verified only for special cases,
\cite{ar:BlaeserDell07,ar:BlaeserDellMakowsky08,ar:BlaeserHoffmann08}.
\item
$Q(G;x_0,y_0)$ is fixed parameter tractable in the sense of
\cite{bk:DowneyFellows99}  when restricted to
graphs classes of bounded tree-width (Proposition \ref{prop:tw})
or even to classes of bounded clique-width (Proposition \ref{prop:cw}).
For the Tutte polynomial, this is known only for graph classes
of bounded tree-width, \cite{ar:Noble98,ar:Andrzejak97,pr:MRAG06}.
\end{itemize}

\subsection*{Outline of the paper}
The paper is organized as follows:
In Section 
\ref{se:Q} we introduce the polynomial $Q(G;x,y)$ and its
univariate versions.
In Seection \ref{se:distinct} we discuss the distinguishing power of
$Q(G;x,y)$ and compare this to other graph polynomials.
In Section \ref{se:combint} we show how certain graph parameters
are definable using $Q(G;x,y)$ and relate it to partition functions
and counting weighted homomorphisms.
In Section \ref{se:decomp} we give a recursive defintion of
$Q(G;x,y)$ using deletion, contraction and extraction of vertices
and show that
$Q(G;x,y)$ is universal. We also compare it to
the universal edge elimination polynomial $\xi(G;x,y,z)$ defined
in \cite{pr:AverbouchGodlinMakowsky07,ar:AverbouchGodlinMakowsky08}
and give a subset expansion formula for
$Q(G;x,y)$. 
In Section \ref{se:cliquesep} we prove decomposition formulas 
for $Q(G;x,y)$ for clique separators.
In Section \ref{se:reconstruct} we show the reconstructibility of
$Q(G;x,y)$. 
In  Section \ref{se:random} we discuss its use to compute
the residual connectedness reliability.
In Section \ref{se:complexity} we discuss the complexity of computing
$Q(G;x,y)$. 
In Section \ref{se:conclu} we draw conclusions and state open problems.
\section{The Subgraph Component Polynomial $Q(G;x,y)$}
\label{se:Q}
\subsection{The Bivariate Polynomial}
Let $G=\left( V,E\right) $ be a finite undirected graph with $n$ vertices
and let $k\leq n$ be a positive integer. Assume the vertices of $G$ fail
stochastic independently with a given probability $q=1-p$. What is the
probability that a subgraph of $G$ with exactly $k$ components survives? The
solution of this problem leads to the enumeration of vertex induced
subgraphs of $G$ with $k$ components. For a given vertex subset $X\subseteq
V $, let $G\left[ X\right] $ be the \emph{vertex induced subgraph} of $G$
with vertex set $X$ and all edges of $G$ that have both end vertices in $X$.
We denote by $k\left( G\right) $ the number of components of $G$. Let 
$q_{ij}\left( G\right) $ be the number of vertex subsets $X\subseteq V$ with 
$i$ vertices such that $G\left[ X\right] $ has exactly $j$ components:
\[
q_{ij}\left( G\right) =\left\vert \left\{ X\subseteq V:\left\vert
X\right\vert =i\wedge k\left( G\left[ X\right] \right) =j\right\}
\right\vert 
\]
The ordinary generating function for these numbers is the two-variable
polynomial
\[
Q\left( G;x,y\right) =\sum_{i=0}^{n}\sum_{j=0}^{n}q_{ij}\left( G\right)
x^{i}y^{j}. 
\]
We call $Q\left( G;x,y\right) $ the 
\emph{subgraph component polynomial} of $G$. 
Since loops or parallel edges do not contribute to connectedness properties of a
graph, we assume in this paper that all graphs are simple.

\begin{figure}[ht]
\begin{center}
\epsfig{file={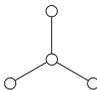},width=0.5in} 
\end{center}
\caption{The star $Star_3=K_{1,3}$} \label{k1_3}
\end{figure}

The star $K_{1,3}$, presented in Figure \ref{k1_3}, has the subgraph
polynomial%
\[
Q\left( K_{1,3};x,y\right)
=1+4xy+3x^{2}y+3x^{3}y+x^{4}y+3x^{2}y^{2}+x^{3}y^{3}. 
\]%
The term $3x^{2}y^{2}$ tell us that there are 3 possibilities to select two
vertices of $G$ that are non-adjacent.

The empty set induces the null graph $N=\left( \emptyset ,\emptyset \right) $
that we consider as being connected, which gives $q_{00}\left( G\right)
=Q\left( G;0,0\right) =1$ for any graph. Substitution of $1$ for $y$ results
in an univariate polynomial that is the ordinary generating function for all
subsets of $V$, i.e. $Q(G;x,1)=\left( 1+x\right) ^{n}$.

\subsection{Univariate Polynomials}
\label{se:univariate}

The coefficient of $y^{k}$ in $Q\left( G;x,y\right) $ is the ordinary
generating function for the number of vertex sets that induce a subgraph of $%
G$ with exactly $k$ components:%
\[
Q_{k}\left( G;x\right) =\left[ y^{k}\right] Q\left( G;x,y\right) 
\]%
We call the polynomial $Q_{k}$ for $k\in \mathbb{N}$ again \emph{subgraph component
polynomial}. The subscript as well as the number of variables should avoid
confusion with the formerly defined subgraph polynomial. The subgraph
polynomial $Q_{1}\left( G;x\right) $ is of special interest. We rename this
polynomial to 
\[
S\left( G;x\right) =Q_{1}\left( G;x\right) =\sum_{i=0}^{n}s_{i}\left(
G\right) x^{i}. 
\]%
It counts the connected vertex induced subgraphs of $G$. A \emph{separating
vertex set} of a connected graph $G=\left( V,E\right) $ is a subset $%
X\subseteq V$ such that $G-X$ is a disconnected graph.

\begin{theorem}
Let $G=\left( V,E\right) $ be a connected graph with $n$ vertices. Let $%
c_{k}\left( G\right) $ be the number of separating vertex sets of
cardinality $k$ for $k=0,1,...,n$. Then the coefficients of the subgraph
polynomial $S\left( G;x\right) $ are given by%
\[
s_{k}\left( G\right) =\dbinom{n}{k}-c_{n-k}\left( G\right) . 
\]
\end{theorem}

\begin{proof}
If $X$ is a separating vertex set then $V\setminus X$ induces a disconnected
graph. Conversely, if $X\subseteq V$ is not a separating vertex set of $G$
then $G\left[ V\setminus X\right] $ is connected.
\end{proof}

We conclude that $2^{n}-S(G;1)$ is the number of all separating vertex sets
of $G$.

\begin{figure}[ht]
\begin{center}
\epsfig{file={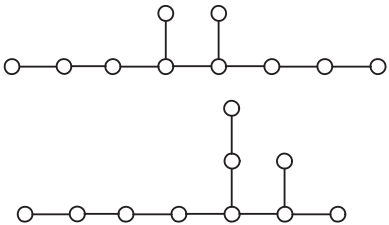},width=2in} 
\end{center}
\caption{Non-isomorphic trees with the same subgraph polynomial} 
\label{tree10}
\end{figure}

A graph invariant $f$ is $trivial$ on a class of graphs $K$
of for any two graphs $G_1$ and $G_2$ with the same number of vertices
we have $f(G_1)= f(G_2)$.

\begin{proposition}
\label{prop:trees}
All non-isomorphic trees with up to nine vertices have different 
subgraph component polynomials. 
In other words we have:
The graph polynomials $Q_k(G;x)$ and $Q(G;x,y)$ are not trivial on trees.
\end{proposition}

However, we have:
\begin{proposition}
\label{prop:trees1}
There exist a unique pair of non-isomorphic trees with
$10$ vertices sharing the same subgraph component polynomial. 
\end{proposition}
\begin{proof}
Figure \ref{tree10} shows
these trees. This statement is true for $S\left( G;x\right) $ as well as for
the (general) subgraph polynomial $Q\left( G;x,y\right) $.
\end{proof}

In Section \ref{se:decomp}
we shall see how to use this to generate infinite families
of pairs of graphs which are not distinguished by $Q(G;x,y)$.

\section{Distinctive Power}
\label{se:distinct}

We denote by $m(G;x)=\sum_i m_i(G) x^i$ 
be the matching polynomial
with $m_i(G)$ the number of $i$-matchings of $G$,
by $p(G;x)$ be the characteristic polynomial, by
$T(G;x,y)$ the Tutte polynomial,
and by
$P(G;x,y)$ the bivariate chromatic polynomial introduced in
\cite{ar:DPT03}.

\begin{figure}[ht]
\begin{center}
\epsfig{file={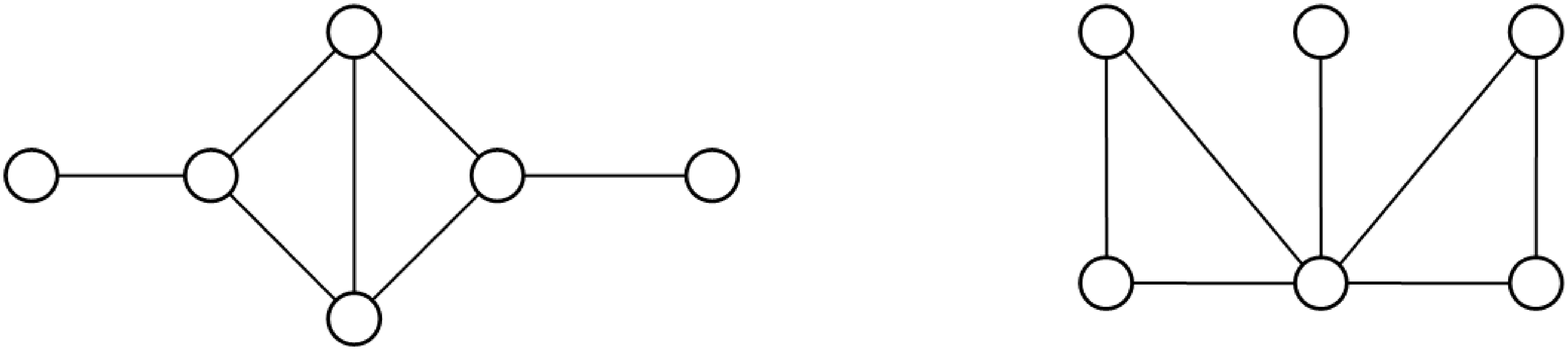},width=3.5in} 
\end{center}
\caption{The graphs $G_1, G_2$} \label{fig:char}
\end{figure}

\begin{figure}[ht]
\begin{center}
\epsfig{file={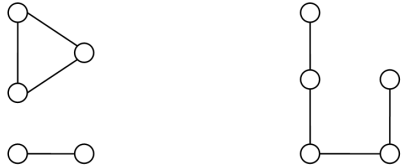},width=2.5in} 
\end{center}
\caption{The graphs $G_3, G_4$} \label{fig:match}
\end{figure}

\begin{figure}[ht]
\begin{center}
\epsfig{file={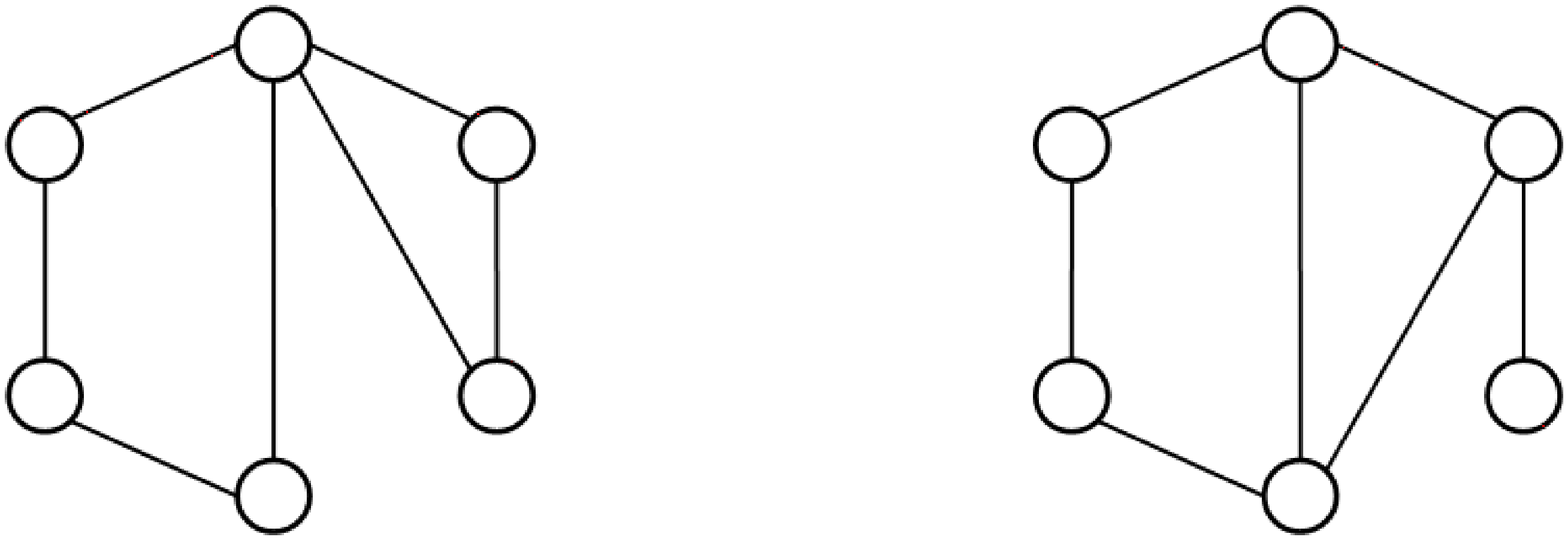},width=2.5in} 
\end{center}
\caption{The graphs $G_5,G_6$} \label{fig:bichrom}
\end{figure}

\begin{proposition}
\begin{enumerate}
For the graphs $G_i; i=1, \ldots 6$ 
from
Figures 
\ref{fig:char}, 
\ref{fig:match}, 
\ref{fig:bichrom},
and for
$P_4$ and $K_{1,3}$ 
we have
\item[(1)]
$p(G_1; x) = p(G_2; x)$
but $Q(G_1;x,y) \neq Q(G_2;x,y)$.
\item[(2)]
$m(G_3; x) = m(G_4; x)$
but $Q(G_3;x,y) \neq Q(G_4;x,y)$.
\item[(3)]
$P(G_5; x,y) = P(G_6; x,y)$
but $Q(G_5;x,y) \neq Q(G_6;x,y)$.
\item[(4)]
$T(P_4; x,y) = T(K_{1,3}; x,y)$
but $Q(P_4;x,y) \neq Q(K_{1,3};x,y)$.
\end{enumerate}
\end{proposition}

\begin{proof}
(1) and (2) are easy to verify.  
\\
For (3) 
$P(G_5; x,y) = P(G_6; x,y)$ is from \cite{ar:DPT03}.
For $Q(G_5;x,y) \neq Q(G_6;x,y)$  we compare 
$[x^4y^3]Q(G_5,x,y)$ with
$[x^4y^3]Q(G_6,x,y)$.
\\
For (4)  we use that the Tutte polynomial does not distinguish
trees of the same size, but that $Q(G;x,y)$ distinguishes all trees
of size up to nine vertices, \ref{prop:trees}.
\end{proof}

\begin{remark}
$Q(G;x,y)$ does not distinguish between graphs which differ only by
the multiplicity of their edges, whereas
for the Tutte polynomial this is not the case.
Let $K^{(m)}_n$ denote the complete graph with $m$ edges between any
two distinct vertices.
Then we have
$T(K^{1}_n;x,y) \neq T(K^{2}_n;x,y)$ but
$Q(K^{1}_n;x,y) = Q(K^{2}_n;x,y)$. 
\end{remark}

\begin{problem}
\label{problem1}
Are there simple graphs
distinguished by $p(G;x)$, $m(G;x)$, $P(G;x,y)$ or $T(G;x,y)$
which are not distinguished by
$Q(G;x,y)$?
\end{problem}

We say that a simple graph $G$ is {\em determined by a graph polynomial} $f$ 
if for every simple graph $G'$
such that $f(G)=f(G')$ we have that $G$ is isomorphic to $G'$.
The class of simple graphs $K$ is {\em determined by a graph polynomial} $f$ 
if every graph $G \in K$ is determined by $f$.
This notion has been studied in
\cite{ar:Noy03,ar:MierNoy04}, 
for the chromatic polynomial, the Tutte polynomial and the
matching polynomial.
It is shown, e.g.,
that several well-known families of graphs are determined by their
Tutte polynomial, among them the class of  wheels, squares of cycles, 
complete multipartite graphs, ladders, M\"obius ladders, and hypercubes.

It follows from Proposition \ref{prop:simple}
that 
the class of empty graphs $E_n$
is determined by
$Q(G;x,y,z)$, and so is the class of complete graphs $K_n$.
Note that, since $T(E_n;x,y)=1$ for all $n \in \mathbb{N}$,
the class of empty graphs is not determined by the
the Tutte polynomial.
It follows from Proposition \ref{prop:trees1}
that 
the class of trees is not determined by $Q(G;x,y,z)$. 

\begin{proposition}
The class of graphs of the form $Star_n= K_{1,n}$ is determined by $Q(G;x,y,z)$. 
\end{proposition}
\begin{proof}
It is easy to verify, that
if 
$[x^{n+1}y]Q(G;x,y) =1$,
$[x^{n}y^{n}]Q(G;x,y) =1$
and
$[x^{n+2}]Q(G;x,y) =0$,
then $G$ is isomorphic to $Star_n$.
\end{proof}

\section{Combinatorial Interpretations}
\label{se:combint}

\subsection{Evaluations and Coefficients of $Q(G;x,y)$}
For a polynomial $f\left( x,y\right) $, let $\left[ x^{i}y^{j}\right]
f\left( x,y\right) $ be the coefficient of $x^{i}y^{j}$ in $f\left(
x,y\right) $ and let $\deg _{x}f$ be the degree with respect to the variable 
$x$. 

\begin{proposition}
\label{prop:simple}
The following graph properties can be easily obtained from the subgraph
polynomial:
\begin{enumerate}
\item[(1)]
The number of vertices:%
\[
n\left( G\right) =\left\vert V\left( G\right) \right\vert =\deg _{x}Q\left(
G;x,y\right) =\log _{2}Q\left( G;1,1\right) 
\]

\item[(2)]
The number of edges:%
\[
e\left( G\right) =\left[ x^{2}y\right] Q\left( G;x,y\right) 
\]

\item[(3)]
The number of components:%
\[
k\left( G\right) =\deg \left( \left[ x^{n\left( G\right) }\right] Q\left(
G;x,y\right) \right) 
\]
\end{enumerate}
\end{proposition}

\begin{theorem}
\label{theo_alpha}The degree of the subgraph component polynomial $Q\left(
G;x,y\right) $ with respect to $y$ is the cardinality of a maximum
independent set of $G$ (the independence number):%
\[
\deg _{y}Q\left( G;x,y\right) =\alpha \left( G\right) 
\]
\end{theorem}

\begin{proof}
Let $X\subseteq V$ be a maximum independent set of $G$. In this case, we
have $k\left( G\left[ X\right] \right) =\left\vert X\right\vert $ and hence $%
\deg _{y}Q\left( G;x,y\right) \geq \alpha \left( G\right) $. Assume that
there exists a set $Y\subseteq V$ with $k\left( G\left[ Y\right] \right)
>\left\vert X\right\vert $. Then we obtain an independent set $X^{\prime }$
by selecting one vertex of each component of $G\left[ Y\right] $ such that $%
\left\vert X^{\prime }\right\vert >\left\vert X\right\vert $ -- a
contradiction.
\end{proof}

Let $a_{i}\left( G\right) $ be the number of independent vertex sets of size 
$i$ of $G$. The \emph{independence polynomial} of $G$ is defined by 
\[
I\left( G;x\right) =\sum_{i=0}^{n}a_{i}\left( G\right) x^{i}. 
\]%
As a consequence of Theorem \ref{theo_alpha} we can derive the independence
polynomial of $G$ from the subgraph component polynomial. Let $\left[ y^{j}\right]
Q\left( G;x,y\right) $ denote the coefficient of $y^{j}$ in $Q\left(
G;x,y\right) $. This coefficient is a polynomial in $x$ where the
coefficient of $x^{i}$ counts the vertex subsets of cardinality $i$ of $G$
that induce a subgraph with $j$ components. A vertex set $X\subseteq V$ is
independent if and only if $k\left( G\left[ X\right] \right) =\left\vert
X\right\vert $. Hence $\left[ x^{j}\right] \left[ y^{j}\right] Q\left(
G;x,y\right) =a_{j}\left( G\right) $ is the number of independent vertex
sets of size $j$ of $G$.

\subsection{Partition Functions}

In this subsection we show that the subgraph polynomial $Q(G;x,y_0)$ 
for any $x\in \mathbb{R}$ and fixed $y_0 \in \mathbb{N}$
can be viewed as a partition function,
using counting of weighted graph homomorphisms. 
Partition functions were first studied in the context of
statistical physics and have recently attracted much attention, \cite{bk:NesetrilWinkler04}.
A systematic study of the question which graph invariants can be
presented as partition functions has been initiated in
\cite{ar:FreedmanLovaszSchrijver07}.

A weighted graph 
$(H, \alpha, \beta)$  
consists of a graph $H=(V(H), E(H))$
with
$\alpha$ assigning weights to vertices and $\beta$ to edges.
The partition function 
$Z_{H, \alpha, \beta}(G)$ 
associated with
$(H, \alpha, \beta)$  
is defined by
$$
Z_H(G) = \sum_{{\tiny
\begin{array}{cc}
h:V\mapsto V_H \\
homomorphism
\end{array}}
}%
\prod_{v\in V}\alpha(h(v))%
\prod_{(u,v)\in E}\beta(h(u),h(v))
$$

\begin{figure}[ht]
\begin{center}
\epsfig{file=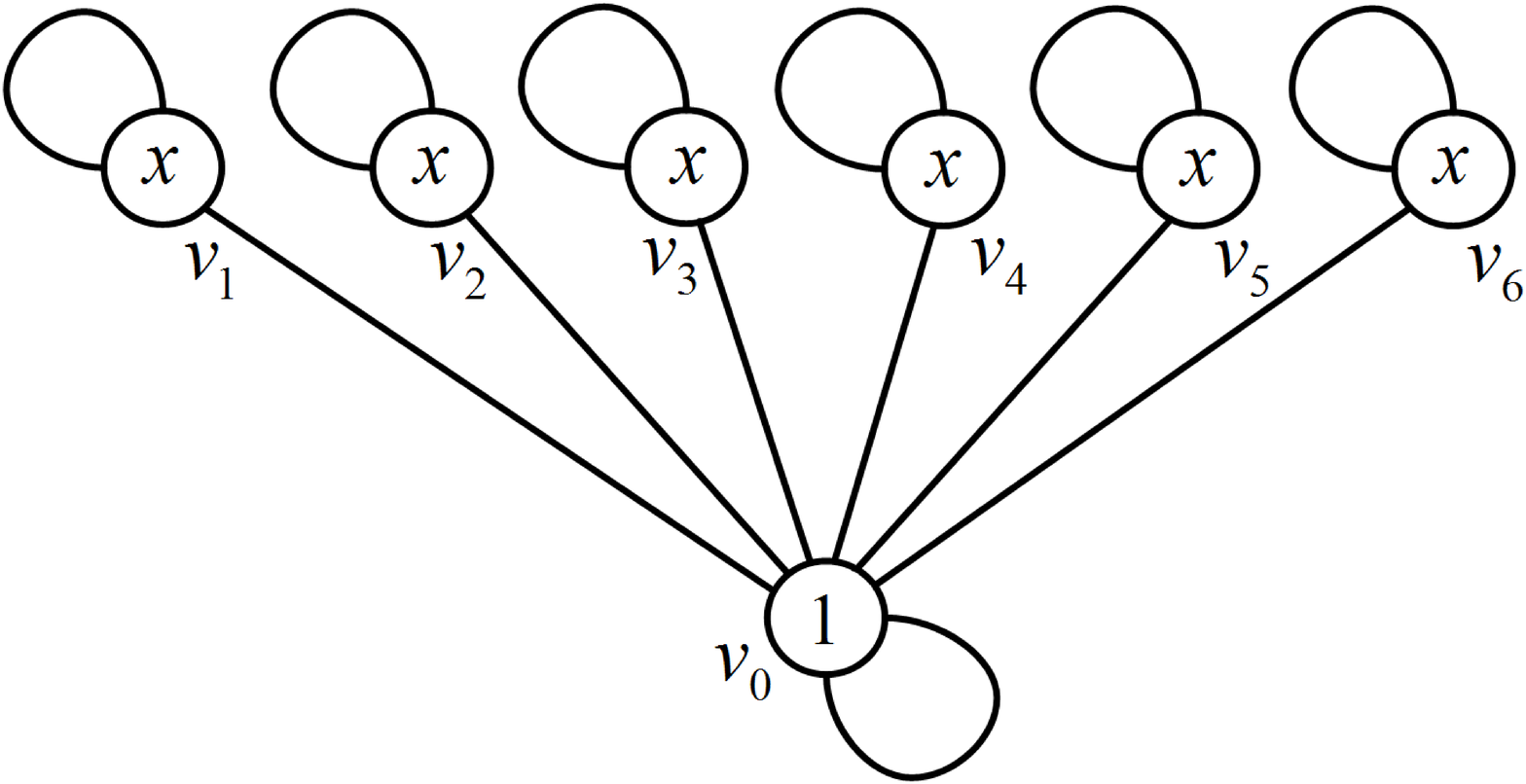,width=0.5\textwidth}
\end{center}
\caption{Auxiliary graph $Star_y$ for $y=6$} \label{fig:hom}
\end{figure}

Let $(Star_y, \alpha, \beta)$ be a weighted star with $y+1$ vertices and with all loops. 
The central vertex is $v_0$.
An example of $Star_y$ for $y=6$ is shown on Figure \ref{fig:hom}.
The weight functions are defined as follows: 

$$ \alpha(v) = \left\{
\begin{array}{ll}
1 & if ~ v=v_0\\
x & otherwise 
\end{array}
\right.$$

$$ \beta(u,v) = 1$$

\begin{theorem}
\label{thm:hom}
Let $Z_{(H, \alpha, \beta)}(G)$ be the partition function associated with
$H = Star_y$ and $\alpha, \beta$ as above.
Then, for all nonnegative integers $y$ and all $x\in \mathbb{R}$, we have
$$
Q(G;x,y) = Z_{(Star_y, \alpha, \beta)}
$$
\end{theorem}
\begin{proof}
Let us start with the definition of $Z_H(G)$. 
Under every mapping $h:V\mapsto V_H$, let 
$A \subseteq V$ be the subset of vertices that are not mapped to $v_0$. 
Let us count the homomorphisms that map the subset $A$ into 
$v_1,\ldots,v_y$: there are exactly $y^{k(G[A])}$ such homomorphisms, because 
every connected component of $G[A]$ must be mapped into a single vertex. 
Finally, we get 
\begin{eqnarray*}
Z_H(G) 
&=& 
\sum_{
{\tiny
\begin{array}{cc}
h:V\mapsto V_H \\
homomorphism
\end{array}
}
}%
\prod_{v\in V}\alpha(h(v)) 
= 
\\
= 
\sum_{{\tiny
\begin{array}{cc}
h:V\mapsto V_H \\
homomorphism
\end{array}}
}%
x^{|A|}  
&=&
\sum_{A\subseteq V} 
\sum_{{\tiny
\begin{array}{c}
h:V\mapsto V_H \\
homomorphism \\
v\in A \leftrightarrow h(v) \neq v_0
\end{array}}
}%
x^{|A|}  
=
\\
= 
\sum_{A\subseteq V} y^{k(G[A])}x^{|A|}
&=&
Q(G;x,y) 
\end{eqnarray*}

which by (\ref{def:subset}) completes the proof
\end{proof}

It is open whether there are other points in which $Q(G;x,y)$ is 
definable as a partition function. 

\begin{remark}
The auxiliary graph $H=Star_n$ is called in physical literature 
{\em The Widom-Rowlinson model for $n$ particles}. The
homomorphisms to $Star_n$ are called {\em Widom-Rowlinson configurations},
\cite{ar:DyerGreenhill2000}.  
\end{remark}

\section{Recursive Definition and Subset Expansion}
\label{se:decomp}

\subsection{Recurrence Relation for Vertex Elimination}

We turn now our attention to the investigation of properties of the subgraph
polynomial that support its computation. The first statement concerns the
multiplicativity with respect to components of the graph.

\begin{theorem}[Multiplicativity]
\label{theo_product}
\begin{enumerate}
\item[(1)]
Let $G = G_1 \sqcup G_1$ be the disjoint union of the graphs $G_1$ and $G_2$.
Then
$$ Q( G;x,y) = Q(G_1;x,y) \cdot Q(G_2;x,y).$$
\item[(2)]
In particular, 
if $G=\left( V,E\right) $ consists of $c$ components 
$G_{1},...,G_{c}$ then the subgraph polynomial satisfies%
\[
Q\left( G;x,y\right) =\prod\limits_{j=1}^{c}Q\left( G_{j};x,y\right) . 
\]
\end{enumerate}
\end{theorem}

\begin{proof}
In case $c=2$, each subset $X\subseteq V$ of cardinality $k$ is the disjoint
union of two subsets $X_{1}\subseteq V\left( G_{1}\right) $ and $%
X_{2}\subseteq V\left( G_{2}\right) $ with $\left\vert X_{1}\right\vert =j$
and $\left\vert X_{2}\right\vert =i-j$. The number of components of $G\left[
X\right] =G\left[ X_{1}\cup X_{2}\right] $ is the sum of the number of
components of $G\left[ X_{1}\right] $ and $G\left[ X_{2}\right] $. We obtain%
\begin{equation}
q_{ik}\left( G\right) =\sum_{j=0}^{i}\sum_{l=0}^{k}q_{jl}\left( G_{1}\right)
q_{i-j,k-l}\left( G_{2}\right) .  \label{product_proof}
\end{equation}%
Thus for a graph with two components, the subgraph polynomial satisfies%
\[
Q\left( G;x,y\right) =Q\left( G_{1};x,y\right) Q\left( G_{2};x,y\right) . 
\]%
We obtain the statement of the theorem by induction on the number of
components.
\end{proof}

\begin{figure}[ht]
\begin{center}
\epsfig{file={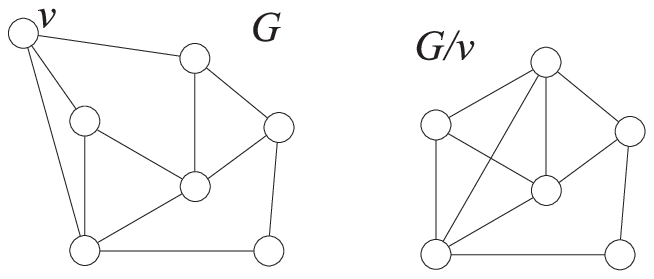},width=2.5962in} 
\end{center}
\caption{Vertex contraction} \label{vertex_contraction}
\end{figure}

We distinguish three types of {\em vertex elimination}:
\begin{description}
\ifcpc
\item[Vertex deletion:]
\else
\item[Vertex deletion]
\fi 
For a given vertex $v\in V\left( G\right) $, let $G-v$ the graph obtained
from $G$ by removal of $v$ and all edges that are incident to $v$.
We call this operation {\em vertex deletion}.
\ifcpc
\item[Vertex extraction:]
\else
\item[Vertex extraction]
\fi 
Similarly, let $G-X$ be the graph obtained from $G$ by removal of all
vertices of the set $X\subseteq V$. Let $N\left( v\right) $ be the set of
vertices that are adjacent to $v$ in $G$ (the neighborhood of $v$). We
denote by $N\left[ v\right] $ the \emph{closed neighborhood} of a vertex $v$
in $G$, i.e. the set of all vertices adjacent to $v$ including $v$ itself. 
The operation
$G-N[v]$ is called
{\em vertex extraction}.
\ifcpc
\item[Vertex contraction:]
\else
\item[Vertex contraction]
\fi 
A further
special graph operation is needed here -- the \emph{contraction} of a
vertex. That is the graph $G/v$ obtained from $G$ by removal of $v$ and
insertion of edges between all pairs of non-adjacent neighbor vertices of $v$. 
Figure \ref{vertex_contraction} shows an example graph and the graph
obtained by vertex contraction.
\end{description}

\begin{theorem}
\label{theo_decom}
Let $G=\left( V,E\right) $ be a graph and $v\in V$. 
Then
the subgraph polynomial satisfies the decomposition formula
\[
Q\left( G;x,y\right) =Q\left( G-v;x,y\right) +x\left( y-1\right) Q\left( G-N 
\left[ v\right] ;x,y\right) +xQ\left( G/v;x,y\right) . 
\]
\end{theorem}

\begin{proof}
Let us first consider all vertex induced subgraphs of $G$ that do not
contain vertex $v$. These subgraphs are also vertex induced subgraphs of $%
G-v $. Consequently, 
\[
Q\left( G-v;x,y\right) 
\]
enumerates all induced subgraphs not including the vertex $v$.

In a second step we count all vertex induced subgraphs that contain vertex $%
v $ but none of its neighbors in $G$. In this case, the vertex $v$ forms a
connected component consisting of $v$ only. The rest of the induced subgraph
is a subgraph of of $G-N\left[ v\right] $. All these subgraphs are
enumerated by $Q\left( G-N\left[ v\right] ;x,y\right) $. However, the
component built by $v$ contributes one vertex and one component to the
polynomial. Thus we obtain the generating function 
\[
xyQ\left( G-N\left[ v\right] ;x,y\right) . 
\]

In our enumeration so far we missed exactly those vertex induced subgraphs
that contain $v$ and at least one of its neighbors together in one
component. We include $v$ in the corresponding candidate set, remove it from 
$G$, and multiply the generating function by $x$ (not by $xy$ because we do
not increase the number of components). In order to trace the components, we
have to simulate the paths using $v$ in $G$. These paths are no longer
present in $G-v$. This task is best performed by using $G/v$ instead of $G-v$%
. Thus we obtain the contribution $xQ\left( G/v;x,y\right) $ to the
generating function. Unfortunately, this polynomial enumerates induced
subgraphs that do not contain any vertices from $N\left( v\right) $, too. We
can fix this problem by subtraction of $xQ\left( G-N\left[ v\right]
;x,y\right) $, which gives%
\[
xQ\left( G/v;x,y\right) -xQ\left( G-N\left[ v\right] ;x,y\right) 
\]%
as final contribution to the generating function.
\end{proof}

\begin{corollary}
Let $v\in V$ be a vertex of degree 1 in $G=\left( V,E\right) $ and let $w$
be its only neighbor in $G$. Then%
\[
Q\left( G;x,y\right) =\left( 1+x\right) Q\left( G-v;x,y\right) +x\left(
y-1\right) Q\left( G-\left\{ v,w\right\} ;x,y\right) . 
\]
\end{corollary}

\begin{proof}
Notice that in this case $G/v=G-v$ and $G-N[v]=G-\left\{ v,w\right\} $. Then
the statement follows immediately from Theorem \ref{theo_decom}.
\end{proof}

\subsection{Some Easy Computations}
The subgraph polynomial can be easily computed for certain special graphs.

\begin{proposition}
\label{prop:easy-1}
For the complete graphs $K_n$ and the empty graphs
$E_n$
(the complement of $K_{n}$) 
we have:
\begin{enumerate}
\item[(1)]
$ Q\left( K_{n};x,y\right) =y\left( 1+x\right) ^{n}-y+1$.
\item[(2)]
$ Q\left( E_{n};x,y\right) =\left( 1+xy\right) ^{n}$.
\end{enumerate}
\end{proposition}
\begin{proof}
(1)
In a complete graph $K_{n}$ each vertex subset except the empty set induces
a connected subgraph. 

(2)
In the empty graph $E_n$ 
each subset $X\subseteq V$
induces a subgraph with $\left\vert X\right\vert $ components.
\end{proof}

From Theorem \ref{theo_decom} we obtain a recurrence relation for the
subgraph polynomial of the paths $P_n$:
\begin{proposition}
\label{prop:easy-2}
\[
Q\left( P_{n};x,y\right) =\left( 1+x\right) Q\left( P_{n-1};x,y\right)
+x\left( y-1\right) Q\left( P_{n-2};x,y\right)  \label{path}
\]
\end{proposition}
Together with the initial values,
\begin{eqnarray*}
Q\left( P_{0};x,y\right) &=&1, \\
Q\left( P_{1};x,y\right) &=&1+xy,
\end{eqnarray*}
equation (\ref{path}) determines the subgraph polynomial of $P_{n}$
uniquely. The explicit solution is
\[
Q\left( P_{n};x,y\right) =\frac{1-x+a}{2a}\left( \frac{2x\left( 1-y\right) }{%
1+x-a}\right) ^{n+1}-\frac{1-x-a}{2a}\left( \frac{2x\left( 1-y\right) }{1+x+a%
}\right) ^{n+1} 
\]%
with $a=\sqrt{1-2x+x^{2}+4xy}$.

The subgraph polynomial of the cycle $C_{n}$ satisfies another recurrence
relation:

\begin{proposition}
\label{prop:easy-3}
\[
Q\left( C_{n};x,y\right) =Q\left( P_{n-1};x,y\right) +x\left( y-1\right)
Q\left( P_{n-3};,x,y\right) +xQ\left( C_{n-1};x,y\right) 
\]
\end{proposition}

The \emph{join} $G\vee H$ of two graphs $G=\left( V,E\right) $ and $H\left(
W,F\right) $ with $V\cap W=\emptyset $ is the graph obtained from $G\cup H$
by introducing edges form each vertex of $G$ to each vertex of $H$.
Consequently, the join of two empty graphs $\overline{K_{s}}$ and $\overline{%
K_{t}}$ is the complete bipartite graph $K_{s,t}$.

\begin{theorem}
\label{th:join}
\label{theo_join}Let $G=\left( V,E\right) $ and $H\left( W,F\right) $ be two
graphs with $V\cap W=\emptyset $, $\left\vert V\right\vert =s$, $\left\vert
W\right\vert =t$. Then%
\[
Q\left( G\vee H;x,y\right) =Q\left( G;x,y\right) +Q\left( H;x,y\right) + 
\left[ \left( 1+x\right) ^{s}-1\right] \left[ \left( 1+x\right) ^{t}-1\right]
y-1. 
\]
\end{theorem}

\begin{proof}
All vertex subsets of $V\cup W$ belong to exactly one of three classes:
\newline
(1) subsets of $V$,\newline
(2) subsets of $W,$\newline
(3) subsets that have at least one vertex of $V$ and at least one vertex of 
$W$.

The first two classes are counted by $Q\left( G;x,y\right)$ and 
$Q\left( H;x,y\right)$, respectively. The empty set is counted twice, which is
corrected by subtracting one. All vertex subsets of the third class induce
connected subgraphs of $G\vee H$. The generating function for the number of
subsets of this class is $\left[ \left( 1+x\right) ^{s}-1\right] \left[
\left( 1+x\right) ^{t}-1\right] $.
\end{proof}

From Theorem \ref{theo_join}, we deduce the subgraph polynomial of the
complete bipartite graph:
\begin{corollary}
\[
Q\left( K_{s,t};x,y\right) =\left( 1+xy\right) ^{s}+\left( 1+xy\right) ^{t}+ 
\left[ \left( 1+x\right) ^{s}-1\right] \left[ \left( 1+x\right) ^{t}-1\right]
y-1 
\]
\end{corollary}

Propositions 
\ref{prop:easy-2}
and
\ref{prop:easy-3}
and Theorem \ref{theo_join}
are explicit instances of general results which follow from
the fact that $Q(G;x,y)$ is definable in Monadic Second Order Logic.
We shall discuss this feature in Section \ref{se:msol}.

We can use Theorem \ref{th:join} and the multiplicativity of $Q(G;x,y)$
to prove the following:
\begin{proposition}
\label{prop:inftrees}
There are infinite families 
of pairs of non-isomorphic graphs with a fixed number of
connected components which are not distinguished by $Q(G;x,y)$.
\end{proposition}
\begin{proof}
Let $G$ be a graph.
We define inductively
\begin{eqnarray*}
F_0(G) &=& G\\
J_0(G) &=& G \\
F_{n+1}(G) &=&F_n(G) \sqcup G\\
J_{n+1}(G) &=&J_n(G) \vee G
\end{eqnarray*}

Let $Tr_1$ and $Tr_2$ be the two trees from Figure \ref{tree10}.
Then,
using
Theorem \ref{th:join} and the multiplicativity of $Q(G;x,y)$ we have 
for all $n \in \mathbb{N}$
$$
Q(F_n(Tr_1);x,y)=
Q(F_n(Tr_2);x,y)
$$
and
$$
Q(J_n(Tr_1);x,y)=
Q(J_n(Tr_2);x,y)
$$
For $G$ connected,
the graphs 
$J_{n}(G)$ are connected.
The graphs
$F_{n}(G)$ have exactly $n$ components.
So for $m$ components we combine 
$F_{m-1}(G) \sqcup J_n(G)$  which has $m$ components.
\end{proof}

\subsection{The Universality Property of $Q\left( G;x,y\right)$}
\label{se:uniprop}

The vertex decomposition formula represented in Theorem \ref{theo_decom} can
be considered as a vertex equivalent to the well-known edge decomposition
(deletion-contraction relations). Edge decomposition formulae of the form $%
f\left( G\right) =\alpha \left( e\right) f\left( G-e\right) +\beta \left(
e\right) f\left( G/e\right) $ apply to the Tutte polynomial and derived
graph invariants, for instance the number of spanning trees or the
reliability polynomial. Indeed, it was shown by J.G. Oxley and D.J.A. Welsh,
\cite{ar:OxleyWelsh79},  that the Tutte polynomial is
in a certain sense {\em universal}, meaning that all other graph invariants that
satisfy edge decomposition formulae can be derived from the Tutte polynomial
by substitution of variables. A textbook presentation is given in 
\cite{bk:Bollobas99}.
A general framework analyzing universality properties of graph polynomials
is studied in \cite{pr:GodlinKatzMakowsky07}.

It seems natural to ask for the most general
vertex decomposition formula. Let us assume that we try to construct an
ordinary generating function $f\left( G\right) $ that counts some type of
vertex induced subgraphs with respect to the number of vertices. Which
properties should such a function have? If the subgraphs in question are
composed from subgraphs of the components then we can expect multiplicativity of 
$f$ with respect to components of the graph. In order to assign the value 
$f\left( G\right)$ uniquely to a graph $G$ by application of a decomposition
formula as given in Theorem \ref{theo_decom}, certain initial values for the
null graph and the empty graph have to be given. Therefore, we presuppose
the following four properties of $f$:

\begin{enumerate}
\item[(a)] (Multiplicativity) If $G_{1}$ and $G_{2}$ are components of $G$ then $f\left(
G\right) =f\left( G_{1}\right) f\left( G_{2}\right) $.

\item[(b)] (Recurrence relation) Let $\alpha ,\beta ,\gamma \in \mathbb{R}$ and let $v$ be a
vertex of $G$, then%
\begin{equation}
f\left( G\right) =\alpha f\left( G-v\right) +\beta f\left( G-N\left[ v\right]
\right) +\gamma f\left( G/v\right) .  \label{ansatz}
\end{equation}

\item[(c)] (Initial condition) There exists $\delta \in \mathbb{R}$ such that $f\left( \emptyset
\right) =\delta $ for the null graph $\emptyset =\left( \emptyset ,\emptyset
\right) $.

\item[(d)] (Initial condition) There exists $\varepsilon \in \mathbb{R}$ such that $f\left(
E_{1}\right) =\varepsilon $ for a graph $E_{1}=\left( \left\{ v\right\}
,\emptyset \right) $ consisting of one vertex.
\end{enumerate}

Furthermore, in order to make $f$ a well-defined graph polynomial, the result
of computing $f$ has to be the same, irrespective of the order
in which we apply the enabled computation steps.
In particular, it has to be {\em independent of the order of the vertices},
which we use to apply the decomposition formula (b).
In general we may choose $\alpha ,\beta ,\gamma ,\delta ,\varepsilon $ from
a field of characteristic zero or from a ring. 
A graph invariant is {\em proper} if there are two graphs $G_1$ and $G_2$
with the same number of vertices such that $f(G_1) \neq f(G_2)$.

Applying the conditions (b),
(c), and (d) we obtain from $E_{1}-v=E_{1}-N\left[ v\right]
=E_{1}/v=\emptyset $ the equation%
\[
\varepsilon =\left( \alpha +\beta +\gamma \right) \delta . 
\]%
Computing $f\left( E_{2}\right) =f\left( \overline{K_{2}}\right) $ in two
ways using (a) and (b), respectively, results in%
\[
\varepsilon ^{2}=\left( \alpha +\beta +\gamma \right) \varepsilon . 
\]%
Consequently, the values of $\delta $ and $\varepsilon $ are determined:%
\begin{eqnarray*}
\delta &=&1 \\
\varepsilon &=&\alpha +\beta +\gamma
\end{eqnarray*}%
If the constants $\alpha ,\beta ,\gamma $ are properly defined then the
value of $f\left( G\right) $ does not depend on the choice of the vertex $v$
in equation (\ref{ansatz}). Consequently, the function $f\left( G\right) $
does not depend on the order of the vertex decomposition (\ref{ansatz}). The
calculation of $f\left( P_{3}\right) $ for path of three vertices yields, in
case we start from a vertex of degree 1,%
\[
f\left( P_{3}\right) =\left( \alpha +\gamma \right) ^{2}\left( \alpha +\beta
+\gamma \right) +\beta \left( \alpha +\gamma \right) +\beta \left( \alpha
+\beta +\gamma \right) . 
\]%
If we begin the vertex decomposition at the vertex of degree 2 then we obtain%
\[
f\left( P_{3}\right) =\alpha \left( \alpha +\beta +\gamma \right) ^{2}+\beta
+\gamma \left( \alpha +\gamma \right) \left( \alpha +\beta +\gamma \right)
+\beta \gamma 
\]%
These two results coincide if $\beta =0$, $\alpha =1$, or $\alpha +\beta
+\gamma =1$. In any case, there remain only two variables that can be chosen
independently. In case of $\beta =0$, all graphs with the same number of
vertices result in the same polynomial. Therefore, this case does not yield
any interesting applications. If $\alpha +\beta +\gamma =1$ then $f\left(
G\right) =1$ for all graphs. The only remaining choice, $\alpha =1$, gives
for $\beta =x\left( y-1\right) $ and $\gamma =x$ the subgraph component polynomial
$Q(G;x,y)$.

From this we get that 
$Q(G;x,y)$ is universal among polynomials recursively defined
using vertex deletion, vertex extraction and vertex contraction.
More precisely, we have the following theorem.
\begin{theorem}[Universality of $Q(G;x,y)$]
\label{th:universal}
\begin{enumerate}
\item[(1)]
For a graph polynomial $f(G; \alpha, \beta, \gamma, \delta, \varepsilon)$ to be proper and
well-defined by conditions (a)-(d) we have
$\alpha =1$, $\delta=1$ and $ \varepsilon= 1 + \beta +\gamma$.
\item[(2)]
There is a unique proper graph polynomial $U(G; \beta, \gamma)$
which is well-defined by conditions (a)-(d)
and we have
$$
Q(G; x,y ) = U(G; x(y-1), x)
$$
and
$$
U(G; \beta, \gamma) = Q(G; \gamma, \frac{\beta}{\gamma}+1)
$$ 
\end{enumerate}
\end{theorem}

\subsection {Vertex Eliminations vs Edge Elimination}
\ifskip
\else
The subgraph component polynomial $Q(G;x,y)$ can be regarded as counting 
vertex set expansions.  
In the literature there is a variety of graph polynomials, 
including the Tutte polynomial, which can be 
defined by counting edge set expansions. 
\fi 

We have seen in Theorem \ref{th:universal} that
$Q(G;x,y)$ is 
universal among the polynomials
defined recursively via deletion, extraction and contraction of vertices.
In 
\cite{pr:AverbouchGodlinMakowsky07,ar:AverbouchGodlinMakowsky08}
the polynomial $\xi(G;x,y,z)$ was shown to be 
universal among the polynomials
defined recursively via deletion, extraction and contraction of edges.
In this section we will show the connection of $G(G;x,y)$ to 
the universal
edge elimination polynomial $\xi(G;x,y,z)$. 

The polynomial $\xi(G;x,y,z)$ generalizes both the Tutte and the matching polynomials, 
as well as the bivariate chromatic polynomial of \cite{ar:DPT03}.
We shall use the recursive decomposition of $\xi(G;x,y,z)$ 
from \cite{ar:AverbouchGodlinMakowsky08}:

\begin{eqnarray}
\nonumber &&\xi(G;x,y,z) = \xi(G-e;x,y,z) + y\xi(G/e;x,y,z) + z\xi(G\dagger e;x,y,z)\\
\nonumber &&\xi(G_1 \sqcup G_2;x,y,z) = \xi(G_1;x,y,z)\xi(G_2;x,y,z) \\
\nonumber &&\xi(E_1;x,y,z) = x \\
&&\xi(\emptyset) = 1
\label{def:rec_xi}
\end{eqnarray}
where $G_1 \sqcup G_2$ denotes the disjoint union of graphs $G_1$ and $G_2$, 
and the three edge elimination operations are defined as follows: 
\begin{description}
\item[Edge deletion:] We denote by $G-e$ the graph obtained from $G$ by
simply removing the edge $e$.
\item[Edge extraction:] We denote by $G \dagger e$ the graph induced by
$V\setminus\{u,v\}$ provided $e=\{u,v\}$. Note that this operation
removes also all the edges adjacent to $e$.
\item[Edge contraction:] We denote by $G/e$ the graph obtained
from $G$ by unifying the endpoints of $e$.
\end{description}

We will rewrite the decomposition of $Q(G;x,y)$ using Theorem \ref{theo_decom}.

\begin{eqnarray}
\nonumber &&Q(G;x,y) =  Q(G-v;x,y) + x Q(G/v;x,y) + x(y-1)  Q(G-N[v];x,y)\\
\nonumber &&Q(G_1 \sqcup G_2; x,y) = Q(G_1;x,y) Q(G_2;x,y) \\
\nonumber &&Q(E_1;x,y) = xy+1 \\
&&Q(\emptyset) = 1
\label{def:rec_q}
\end{eqnarray}
\begin{theorem}{\label{thm:q_xi}}
Let $G=(V,E)$ be a graph. Let $L(G)=(V_e,E_e)$ denote the line graph of $G$.
Then the following equation holds: 
$$\xi(G;1,x,x(y-1)) = Q(L(G); x, y)$$
\end{theorem}
\begin{proof}
First, let us analyze the correspondence of the edge elimination operations in a graph 
to the vertex elimination operations in its line graph. 
Let $v_e\in V_e$ be the vertex in the line graph that corresponds to the edge $e\in E$
of the original graph. By the definition of the edge and vertex elimination operations: 
\begin{eqnarray}
\label{line_1}&& L(G-e) = L(G)-v_e \\
\label{line_2}&& L(G/e) = L(G)/v_e \\
\label{line_3}&& L(G\dagger e) = L(G)-N[v_e]
\end{eqnarray}
First let us check the connected graphs with up to one edge: 

If $G \in \{\emptyset, E_1\}$, $L(G)=\emptyset$, \\
The equivalence $\xi(G;1,x,x(y-1)) = 1 = Q(\emptyset)$ holds. 

If $G$ is a single point with a loop, or $G=P_2$, $L(G)$ is a singleton, 
The equivalence $\xi(G;1,x,x(y-1)) = 1 + x + x(y-1) = 1 + xy = Q(E_1)$ holds.

Next, we note that $L(G_1 \sqcup G_2) = L(G_1) \sqcup L(G_2)$. Therefore, if 
the theorem holds for graphs $G_1$ and $G_2$, then it holds also for $G_1 \sqcup G_2$. 
Finally, the theorem follows by induction on the number of edges, 
using the decomposition formulae (\ref{def:rec_q}) and (\ref{def:rec_xi}) and the 
correspondence of edge and vertex elimination operations. 
\end{proof}

\begin{problem}
\label{problem2}
How does the distinguishing power of 
$\xi(G;x,y,z)$ compare to the distinguishing power
of $Q(G;x,y)$?
\end{problem}

\subsection{Subset Expansion and Definability in Logic}
\label{se:msol}

$Q(G;x,y)$ was defined as a generating function.
Let us rewrite the definition of 
$Q(G;x,y)$ 
in a slightly different way.
Instead of summation over the number of the used vertices $i$ , and the number of 
induced connected components $j$, we shall summate over all the possible subsets
of vertices:

\begin{equation}\label{def:subset}
Q(G; x,y) = \sum_{A \subseteq V}x^{|A|}y^{k(G[A])}.
\end{equation}

This is a {\em subset expansion formula}, a term coined in \cite{ar:Traldi04}.
The relationship between recursive definitions of graph polynomials
and the existence of subset expansion formulas has been studied
from a logical point of view in \cite{pr:GodlinKatzMakowsky07}.
Subset expansion formulas can often be used to show that
a graph polynomial is definable in Monadic Second Order Logic, as studied
in 
\cite{ar:MakowskyTARSKI,ar:MakowskyZoo} .
However, the exponent $k(G[A])$ in Equation (\ref{def:subset})
causes a problem.
to remedy this, we use, like in \cite{ar:Makowsky01},
an auxiliary order $\prec$ over the vertices. We will denote
by $F(A)$ the subset of the \textit{smallest} vertices in every respective connected 
component. 

\begin{equation}\label{def:msol}
Q(G; x,y) = \sum_{A \subseteq V}\left(\prod_{v\in A}x\right) \left(\prod_{u\in F(A)}y\right)
\end{equation}
Note that the result does not depend on the used auxiliary order.

Without having to go in the details of graph polynomials definable in
Monadic Second Order Logic\footnote{
The interested reader can consult \cite{bk:FMT} for the use
Monadic Second Order Logic in finite model theory,
and \cite{bk:Courcelle09} for its use in graph theory.
},
Equation (\ref{def:msol}) shows that
$Q(G; x,y)$ is a graph polynomial definable in Monadic Second Order Logic
for graphs $G=(V,E)$ with universe $V$ and a binary edge relation.
Therefore all the theorems from 
\cite{ar:MakowskyTARSKI,ar:FischerMakowsky08} can be applied.
In particular, the Feferman-Vaught-type theorems from
\cite{ar:MakowskyTARSKI}
guarantee existence of reduction formulas
like multiplicativity from Theorem \ref{theo_product}, or
the one in Theorem \ref{th:join} for the join, not only for the disjoint union
or the join operation, but for a wide class of $\MSOL$-definable operations.
Also, a general theorem from
\cite{ar:FischerMakowsky08} guarantees the existence of recurrence formulas,
as proven in Propositions 
\ref{prop:easy-2}
and
\ref{prop:easy-3},
for a wide class of recursively defined families of graphs,
as studied also in \cite{ar:NoyRibo04}.
Among these we have the wheels $W_n$, the ladders $L_n$ and the stars $Star_n$. 
It should not be difficult to compute the recurrence relations for these
explicitly.

We shall exploit $\MSOL$-definability also
for our complexity analysis in Section \ref{se:fpt}.

\section{Clique Separators}
\label{se:cliquesep}

The simplest case of a clique separator in a graph $G$ is an \emph{%
articulation}, i.e. a vertex whose removal from $G$ results in an increase
of the number of components.
Let $v$ be an articulation of $G$ and let $H$ and $K$ be
subgraphs of $G$ such that $G=H\cup K$ and 
$H\cap K=\left( \left\{ v\right\},\emptyset \right) $. 
It is well known, cf. \cite{bk:Bollobas99},
that in this case the Tutte polynomial $T(G;x,y)$ satisfies
$$
T(G;x,y)= T(H;x,y) \cdot T(H;x,y)
$$
In the case of the subgraph component polynomial the situation is a bit more
complicated:
\begin{theorem}
\label{theo_arti}
Let $v$ be an articulation of $G$ and let $H$ and $K$ be
subgraphs of $G$ such that $G=H\cup K$ and $H\cap K=\left( \left\{ v\right\}
,\emptyset \right) $. Then the subgraph polynomial $Q\left( G\right)
=Q\left( G;x,y\right) $ satisfies%
\begin{eqnarray*}
Q\left( G\right) &=&Q\left( H-v\right) Q\left( K-v\right) \\
&&+\frac{1}{xy}\left[ Q\left( H\right) -Q\left( H-v\right) \right] \left[
Q\left( K\right) -Q\left( K-v\right) \right] .
\end{eqnarray*}
\end{theorem}

\begin{proof}
The first product of the polynomial is the generating function for the
number all vertex induced subgraphs that do not contain the articulation $v$. 
The product is justified by Theorem \ref{theo_product}. The second term
counts all remaining subgraphs, i.e. those ones containing vertex $v$. 
Here the equation (\ref{product_proof}) from the proof of Theorem 
\ref{theo_product} has to be modified. The vertex $v$ is counted twice because
it belongs to both $K$ and $H$. This double counting is corrected by
multiplication with $x^{-1}$. By analogy, we introduce the factor $y^{-1}$
in order to avoid that the component containing $v$ is counted twice.
\end{proof}

Theorem \ref{theo_arti} can be generalized in order to cover clique
separators with more than one vertex. Let $G=\left( V,E\right) $ be a
connected graph and $H$, $K$ subgraphs of $G$ such that $H\cap K=K_{r}$ and $%
H\cup K=G$. In this case $K_{r}=\left( U,F\right) $ forms a \emph{separating
clique} of $G$. Here we assume that neither $H$ nor $K$ coincides with $%
K_{r} $. The subgraphs $H$ and $K$ are called \emph{split components} of $G$
with respect to $K_{r}$.

\begin{theorem}
\label{th:splitting}
Let $K_{r}=\left( U,F\right) $ be a clique separator of $G$ such that there
are two split components $H$ and $K$. Then the subgraph polynomial $Q\left(
G\right) =Q\left( G;x,y\right) $ satisfies%
\begin{eqnarray*}
Q\left( G\right) &=&Q\left( H-U\right) Q\left( K-U\right) \\
&&+\frac{1}{y}\sum_{\emptyset \neq A\subseteq U}\frac{1}{x^{\left\vert
A\right\vert }}\sum_{B\supseteq U\setminus A}\sum_{C\supseteq U\setminus
A}\left( -1\right) ^{\left\vert B\right\vert +\left\vert C\right\vert
}Q\left( H-B\right) Q\left( K-C\right) .
\end{eqnarray*}
\end{theorem}

\begin{proof}
First we count all subgraphs that are induced by vertex subsets included in $%
V\setminus U$. These subgraphs are also subgraphs of $G-U$. From Theorem \ref%
{theo_product} we obtain $Q\left( H-U\right) Q\left( K-U\right) $ as
generating function for all subgraphs of $G$ induced by subsets of $%
V\setminus U$.

For each subset $A\subseteq U$, let $f_{ij}\left( H,A\right) $ be the number
of vertex subsets $X\subseteq V\left( H\right) $ of cardinality $i$ with $%
A\subseteq X$ such that the induced subgraph $H\left[ X\right] $ has exactly 
$j$ components:%
\[
f_{ij}\left( H,A\right) =\left\vert \left\{ X:A\subseteq X\subseteq V\left(
H\right) \wedge \left\vert X\right\vert =i\wedge k\left( H\left[ X\right]
\right) =j\right\} \right\vert 
\]%
The polynomial%
\[
F\left( H,A\right) =\sum_{i=0}^{n}\sum_{j=0}^{n}f_{ij}\left( H,A\right)
x^{i}y^{j} 
\]%
is the ordinary generating function for the numbers $f_{ij}\left( H,A\right) 
$. We define the numbers $f_{ij}\left( K,A\right) $ and the corresponding
generating function $F\left( K,A\right) $ for the second split component
analogously. Let $X\subseteq V\left( G\right) $ be a vertex subset with $%
X\cap U=A$. Then the component of $G\left[ X\right] $ that contains $A$ is
counted in $F\left( H,A\right) $ and in $F\left( K,A\right) $. There is
indeed only one component counted twice, since $A$ induces a clique of $H$
and $K$, respectively. Thus we obtain%
\begin{equation}
Q\left( G\right) =Q\left( H-U\right) Q\left( K-U\right) +\frac{1}{y}%
\sum_{\emptyset \neq A\subseteq U}\frac{1}{x^{\left\vert A\right\vert }}%
F\left( H,A\right) F\left( K,A\right) .  \label{eq_Q1}
\end{equation}%
The factor $x^{-\left\vert A\right\vert }$ takes into account that all
vertices of $A$ contribute to $F\left( H,A\right) $ and to $F\left(
K,A\right) $. For each subset $B\subseteq U$, the subgraph polynomial of $%
H-B $ can be represented as a sum of generating functions as follows:%
\[
Q\left( H-B\right) =\sum_{A\subseteq U\setminus B}F\left( H,A\right) 
\]%
We define $\hat{Q}\left( H,U\setminus B\right) =Q\left( H-B\right) $ and
obtain%
\[
\hat{Q}\left( H,U\setminus B\right) =\sum_{A\subseteq U\setminus B}F\left(
H,A\right) 
\]%
or 
\[
\hat{Q}\left( H,B\right) =\sum_{A\subseteq B}F\left( H,A\right) . 
\]%
By M\"{o}bius inversion, we obtain%
\begin{eqnarray*}
F\left( H,A\right) &=&\sum_{B\subseteq A}\left( -1\right) ^{\left\vert
A\right\vert -\left\vert B\right\vert }\hat{Q}\left( H,B\right) \\
&=&\sum_{B\subseteq A}\left( -1\right) ^{\left\vert A\right\vert -\left\vert
B\right\vert }Q\left( H-\left( U\setminus B\right) \right) \\
&=&\sum_{U\setminus B\subseteq A}\left( -1\right) ^{\left\vert A\right\vert
-\left\vert U\setminus B\right\vert }Q\left( H-B\right) \\
&=&\left( -1\right) ^{\left\vert A\right\vert -\left\vert U\right\vert
}\sum_{B\supseteq U\setminus A}\left( -1\right) ^{\left\vert B\right\vert
}Q\left( H-B\right) .
\end{eqnarray*}%
Similarly, we can prove for each $A\subseteq U$ that%
\[
F\left( K,A\right) =\left( -1\right) ^{\left\vert A\right\vert -\left\vert
U\right\vert }\sum_{B\supseteq U\setminus A}\left( -1\right) ^{\left\vert
B\right\vert }Q\left( K-B\right) . 
\]%
The substitution of $F\left( H,A\right) $ and $F\left( K,A\right) $ in
equation (\ref{eq_Q1}) yields%
\begin{eqnarray*}
Q\left( G\right) &=&Q\left( H-U\right) Q\left( K-U\right) \\
&&+\frac{1}{y}\sum_{\emptyset \neq A\subseteq U}\frac{1}{x^{\left\vert
A\right\vert }}\sum_{B\supseteq U\setminus A}\left( -1\right) ^{\left\vert
B\right\vert }Q\left( H-B\right) \sum_{C\supseteq U\setminus A}\left(
-1\right) ^{\left\vert C\right\vert }Q\left( K-C\right) .
\end{eqnarray*}
\end{proof}

\begin{problem}
\label{splitting}
Can we have an analogue of Theorem \ref{th:splitting} for the case
where the separating vertex set is not required to be a clique?
\end{problem}

\section{Reconstruction}
\label{se:reconstruct}

The famous \emph{graph reconstruction conjecture} by Kelly and Ulam \cite{bk:Ulam60} states
the every undirected graph with at least three vertices can be reconstructed
from a deck of its vertex-deleted subgraphs (more precisely from the
corresponding isomorphism classes). See for example the papers 
\cite{ar:BondyHemminger77,ar:Bondy91}
as an introduction into this field. Despite the fact that the
conjecture is still open, many graph invariants and graph polynomials (e.g.
the Tutte polynomial) are known to be reconstructible. We can show that also
the subgraph component polynomial can be reconstructed from the deck of the subgraph
polynomials of the vertex-deleted subgraphs.

\begin{theorem}
\label{th:reconstruct}
The subgraph polynomial $Q\left( G;x,y\right) $ for a graph $G=\left(
V,E\right) $ with $n=\left\vert V\left( G\right) \right\vert \geq 3$ is
uniquely determined by the set of polynomials%
\[
\left\{ Q\left( G-v;x,y\right) :v\in V\left( G\right) \right\} . 
\]%
Let $\hat{\omega}$ be the smallest power of $y$ that appears at least twice
among the terms $x^{n-1}y^{j}$ of the polynomials $Q\left( G-v;x,y\right) $.
Define 
\[
\omega =\left\{ 
\begin{array}{l}
n\text{ if }\hat{\omega}=n-1, \\ 
\hat{\omega}\text{ else.}%
\end{array}%
\right. 
\]%
Then the subgraph polynomial is given by%
\[
Q\left( G;x,y\right) =x^{n}\left[ y^{\omega }+\int_{0}^{\frac{1}{x}%
}t^{n-1}\sum_{v\in V}Q\left( G-v;\frac{1}{t},y\right) dt\right] . 
\]
\end{theorem}

\begin{proof}
Let $k\left( G\right) $ denote the number of components of a graph $G$. In
each graph $G=\left( V,E\right) $ with at least three vertices and at least
one edge there exist two vertices $u,v\in V$ such that $k\left( G-u\right)
=k\left( G-v\right) =k\left( G\right) $. If the term $x^{n-1}y^{j}$ appears
in the polynomial $Q\left( G-v;x,y\right) $ then the number of components of 
$G-v$ equals $j$. Since $k\left( G-v\right) \geq k\left( G\right) $ for each
vertex $v\in V$, the the smallest power of $y$ that appears at least twice
among the terms $x^{n-1}y^{j}$ of the polynomials $Q\left( G-v;x,y\right) $
is equal to $k\left( G\right) $. There is only one exception: If $G$ is the
empty (edgeless) graph then the removal of each vertex of $G$ decreases the
number of components by one, which is taken into consideration within the
definition of $\omega $. Consequently, we obtain $\omega =k\left( G\right) $.

Each vertex induced subgraph with $i<n$ vertices is counted exactly $n-i$
times in the polynomial%
\[
\sum_{v\in V\left( G\right) }Q\left( G-v;x,y\right) . 
\]%
The coefficient of $t^{i}y^{j}$ in%
\[
t^{n-1}\sum_{v\in V}Q\left( G-v;\frac{1}{t},y\right) 
\]%
equals $i$ times the number of vertex induced subgraphs of $G$ with exactly $%
n-i$ vertices and $j$ components. The integration with respect to $t$
transforms $t^{i-1}$ into $\frac{1}{i}t^{i}$ such that the vertex induced
subgraphs are enumerated correctly by the coefficients of the resulting
polynomial. Finally, the bounds of integration and the multiplication with $%
x^{n}$ performs the back-substitution in order to obtain an ordinary
generating function with variables $x$ and $y$.
\end{proof}

\section{Random Subgraphs}
\label{se:random}

Now we assume that the vertices of $G=\left( V,E\right) $ fail stochastic
independently with a given (identical) probability $q=1-p$. We obtain the
probability $P_{k}\left( G\right) $ that a vertex induced subgraph of $G$
has exactly $k$ components from the subgraph polynomial:

\begin{equation}
P_{k}\left( G\right) =\frac{1}{k!}\left. \frac{\partial ^{k}}{\partial y^{k}}%
\left( 1-p\right) ^{n}Q\left( G;\frac{p}{1-p},y\right) \right\vert _{y=0}
\label{prob}
\end{equation}%
The sequence $\left\{ P_{k}\left( G\right) \right\} _{k\in \mathbb{N}}$ is
the distribution of the number of components. Consequently, we obtain%
\[
\sum_{k=0}^{n}P_{k}\left( G\right) =1. 
\]

\begin{figure}[ht]
\begin{center}
\epsfig{file={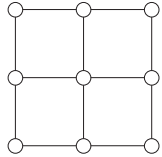},width=0.8in} 
\end{center}
\caption{A $3\times 3$ grid graph} 
\label{grid3_3}
\end{figure}

\begin{figure}[ht]
\begin{center}
\epsfig{file={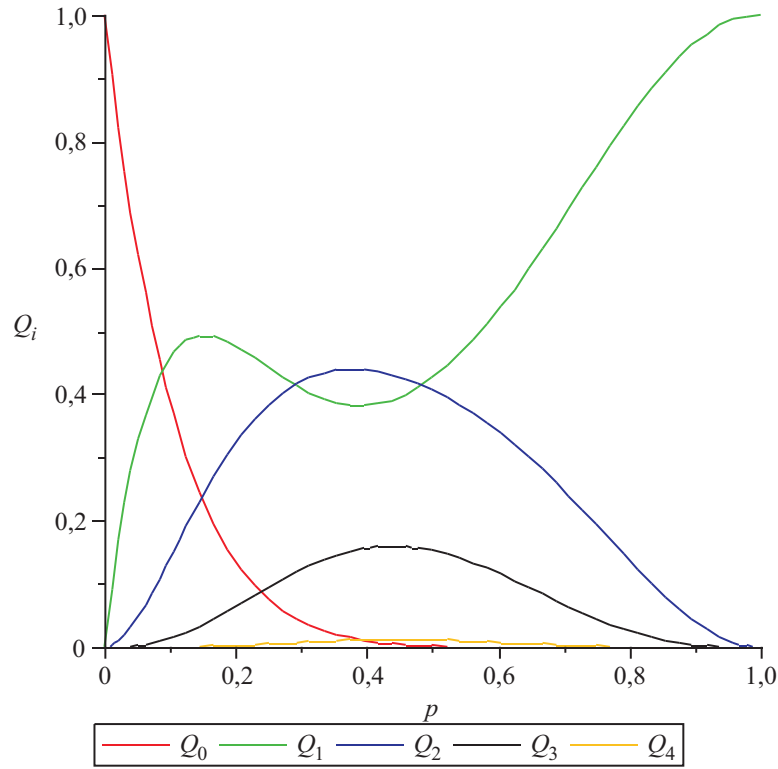},width=4in} 
\end{center}
\caption{Distribution of the number of components} 
\label{probability}
\end{figure}

Figure \ref{probability} shows the distribution for the graph presented in
Figure \ref{grid3_3}.

The probability $P_{1}\left( G\right) $ is called the \emph{%
residual connectedness reliability}. Boesch, Satyanarayana, and Suffel 
\cite{ar:BoeschSatyanarayanaSuffel91}
showed that the computation of $P_{1}\left( G\right) $ is a \#P-hard
problem, even in planar bipartite graphs. Since $P_{1}\left( G\right) $ can
be obtained in polynomial time from the subgraph polynomial by applying the
relation (\ref{prob}), we obtain the following statement.

\begin{corollary}
\label{cor:sharpp}
The computation of the subgraph polynomial is a \#P-hard problem. It remains
\#P-hard for the class of all planar bipartite graphs.
\end{corollary}


\section{Computational complexity of $Q(G;x,y)$}
\label{se:complexity}

\subsection{Complexity of evaluation}
We have already seen in Corollary \ref{cor:sharpp} that $Q(G;x,y)$ 
is $\sharp\mathbf{P}$-hard to compute. 
Now  we deal with a problem of evaluation of $Q(-;x,y)$ at a given point
$(x,y)\in \mathbb{Q}^2$ for arbitrary input graph $G$. 
\begin{theorem}
\label{th:dpp}
For every point $(x,y)\in \mathbb{Q}^2$, possibly except for the lines $xy=0$, 
$y=1$, $x=-1$ and $x=-2$, the evaluation of $Q(G;x,y)$ for an input graph $G$
is  $\sharp\mathbf{P}$-hard. 
\end{theorem}

C. Hoffmann in \cite{pr:Hoffmann08} showed the following:
\begin{theorem}[Hoffmann 2008]
\label{th:hoffmann}
For every point $(x,y,z)\in \mathbb{Q}^3$, except possibly for the subsets
$x=0$, $z=-xy$, $(x,z)\in \{(1,0),(2,0)\}$  and $y\in\{-2,-1,0\}$, 
the evaluation of $\xi(-;x,y,z)$ for an input graph $G$ is $\sharp\mathbf{P}$-hard.
\end{theorem}

\begin{proof}[Proof of Theorem \ref{th:dpp}:]
We use Theorem \ref{th:hoffmann} and our Theorem (\ref{thm:q_xi}). 
Under the conditions of Theorem (\ref{thm:q_xi}), Hoffmann's exception sets are
mapped to the lines  $xy=0$, $y=1$, $x=-1$ and $x=-2$. It follows that for every point
$(x,y)\in \mathbb{Q}^2$ that does not lay on one of those lines, the polynomial
$Q(-;x,y)$ is $\sharp\mathbf{P}$-hard to evaluate even for an input line graph $L(G)$. 
\end{proof}

The evaluation of $Q(-;x,y)$ is polynomial time computable for $xy=0$ and for $y=1$. 
It remains open whether it is polynomial time computable for $x=-1$ and $x=-2$. 
One can also ask, whether there is some point $(x,y)\in \mathbb{Q}^2$, in which 
$Q(-;x,y)$ is hard to evaluate for general input graph, but easy for input line graph. 

\subsection{Parameterized complexity}
\label{se:fpt}
Here we discuss the computational complexity of $Q(G;x,y)$ for input graphs of bounded
tree with, and for input graphs for bounded clique width.
We do not need the exact definitions here.
For background on tree-width the reader can consult
\cite{bk:Diestel}. Clique-width was defined in \cite{ar:CourcelleOlariu00}.
Both are discussed in \cite{ar:MakowskyTARSKI}.

Recall that the subgraph component polynomial is definable using the $\MSOL$-formalism 
(definition \ref{def:msol}) with auxiliary order, while the result is 
order-independent. 
Hense, using a general theorem from \cite{ar:Makowsky01,ar:MakowskyTARSKI},
we have
\begin{proposition}
\label{prop:tw}
$Q(G;x,y)$ is polynomial time computable on graphs
of tree-width at most $k$ where the exponent of the run time is
independent of $k$.
\end{proposition}
Moreover, applying the result of Courcelle, Makowsky and Rotics \cite{ar:CourcelleMakowskyRoticsDAM},
combined with the
results from \cite{ar:Oum2005},
we have a
similar result for graphs of bouded clique width: 
\begin{proposition}
\label{prop:cw}
$Q(G;x,y)$ is polynomial time computable on graphs
of clique-width at most $k$ where the exponent of the run time is
independent of $k$.
\end{proposition}
The drawback of the general methods of
\cite{ar:Makowsky01,ar:MakowskyTARSKI} and \cite{ar:CourcelleMakowskyRoticsDAM},
lies in the huge hidden
constants, which make it practically unusable.
However, an explicit dynamic algorithm for computing the polynomial
$Q(G;x,y)$ on graphs of bounded tree-width, given the tree
decomposition of the graph, where the constants are simply
exponential in $k$, can be constructed along the same ideas as
presented in \cite{ar:Traldi03,ar:FischerMakowskyRavve2006}.
For the graphs of bounded clique width, given the clique decomposition
of the graph, we know an algorithm with constants 
doubly-exponential in $k$. It is open whether an algorithm with constants 
simply exponential in $k$ exists. 
For a comparison of the complexity of computing graph polynomials
on graphs classes of bounded clique-width, cf. \cite{pr:MRAG06}.

\section{Conclusions and Open Problems}
\label{se:conclu}

We have shown that
$Q(G;x,y)$ is a universal vertex elimination polynomial.
We have given a few combinatorial interpretations of its evaluations
and coefficients. 
We have proven various splitting formulas for
$Q(G;x,y)$ such as the multiplicativity, Theorem \ref{theo_decom}
and Theorem \ref{th:splitting}.
Problem \ref{splitting} asks for more such theorems. 
Besides having algorithmic importance, such splitting formulas
increase our structural understanding of the graph polynomial
under study, and may help us in analizing its distictive power.

We have looked at the graph polynomial
$Q(G;x,y)$ 
from
various angles and compared its behaviour and distinguishing
power with the characteristic polynomial, the matching polynomial
the Tutte polynomial and the universal edge elimination polynomial.
We have not discussed
the relationship of
$Q(G;x,y)$ to other graph polynomials, such as the interlace
polynomial, \cite{ar:ArratiaBollobasSorkin2004,ar:AignervdHolst2004},
or the many other graph polynomials listed in \cite{ar:MakowskyZoo}.

We have seen that 
$Q(G;x,y)$ 
distinguishes between graphs
where these polynomials do not.
We have not found cases where these other polynomials do
distinguish between graphs  where 
$Q(G;x,y)$ does not. This is probably due to our 
lack of computerized tools for searching for such cases,
cf. Problem \ref{problem1}.
In Problem \ref{problem2}
we ask about comparing distinguishing power of 
$Q(G;x,y)$  and the universal edge elimination polynomial $\xi(G;x,y,z)$.
This seems to be more tricky.
We have given a few examples of graphs and graph families which are  determined by
$Q(G;x,y)$. 

\begin{problem}
\label{problem3}
Find more graph invariants which are determined by
$Q(G;x,y)$ .
\end{problem}

\begin{problem}
\label{problem4}
Find more classes of graphs which are determined by $Q(G;x,y)$.
\end{problem}

Returning to our motivation, we have only studied the simplest case
of community structure in networks. We have studied the generating function
of induced subgraphs with $i$ vertices which have $j$ components.
More generally, one would want to study community structures
where components are replaced by maximal $k$-connected components.

\begin{problem}
What are the appropriate generating functions which capture the
essence of various community structures?
\end{problem}

\end{document}